\documentclass{amsart}

\usepackage[dvipsnames]{xcolor}
\usepackage[english]{babel}
\usepackage{amsfonts}
\usepackage{amsthm}
\usepackage{amsbsy}
\usepackage{amssymb}
\usepackage{graphicx}
\usepackage{eso-pic}
\usepackage{tikz}
\usepackage{xfrac}
\usepackage{float}
\usepackage[bookmarks=true, colorlinks=true, urlcolor=Blue,linkcolor=BrickRed,citecolor=MidnightBlue,pagebackref=true]{hyperref}\usepackage{bookmark}
\usepackage{amsmath}
\usepackage{subfigure}
\usepackage{hyperref}
\usepackage{scalerel,stackengine}
\usepackage{epigraph}
\usepackage{placeins}
\usetikzlibrary{matrix,arrows}
\usepackage{geometry}
\usepackage{enumitem}

\usepackage{mathtools}
\usepackage{hhline}

\usepackage{amsmath,amssymb,tikz-cd}
\usepackage[inline]{asymptote}

\makeatletter
\newcommand{\mylabel}[2]{#2\def\@currentlabel{#2}\label{#1}}
\makeatother
\usepackage{todonotes}

\newtheorem{lemma}{Lemma}
\newtheorem{thm}[lemma]{Theorem}
\newtheorem{prop}[lemma]{Proposition}
\newtheorem{cor}[lemma]{Corollary}
\newtheorem*{cor*}{Corollary}
\newtheorem*{thm*}{Theorem}

\theoremstyle{definition}

\theoremstyle{remark}

\newcommand\restr[2]{{
  \left.\kern-\nulldelimiterspace 
  #1 
  \littletaller 
  \right|_{#2} 
  }}

\newcommand{\littletaller}{\mathchoice{\vphantom{\big|}}{}{}{}}

\newcommand{\matR}{\mathbb{R}}
\newcommand{\matQ}{\mathbb{Q}}
\newcommand{\matZ}{\mathbb{Z}}
\newcommand{\matC}{\mathbb{C}}

\newcommand{\matH}{\mathbb{H}}

\newcommand{\PO}{\mathrm{PO}}
\newcommand{\PSO}{\mathrm{PSO}}

\newcommand{\Isom}{\mathrm{Isom}}

\newcommand{\thick}{\mathrm{thick}}

\newcommand{\w}[1]{{\color{white} #1}}

\author{Stefano Riolo}
\address{Dipartimento di Matematica, Università di Bologna \newline Piazza di Porta San Donato 5, 40126 Bologna, Italy}
\email{stefano\w{0}.\w{0}riolo\w{0}@unibo.it}
\urladdr{\href{https://www.dm.unibo.it/~stefano.riolo}{www.dm.unibo.it/~stefano.riolo}}

\author{Edoardo Rizzi}
\address{Scuola Normale Superiore \newline Piazza dei Cavalieri 7, 56126 Pisa, Italy}
\email{edoardo\w{0}.\w{0}rizzi\w{0}@sns.it}
\urladdr{\href{https://www.sns.it/en/persona/edoardo-rizzi}{www.sns.it/en/persona/edoardo-rizzi}}

\thanks{\tiny S.R. was funded by the European Union – NextGenerationEU under the National Recovery and Resilience Plan (PNRR) -- Mission 4 Education and research -- Component 2 From research to business -- Investment 1.1 Notice Prin 2022 -- DD N. 104 del 2/2/2022, with title ``Geometry and topology of manifold'', proposal code 2022NMPLT8 -- CUP J53D23003820001.}

\title{A cusped hyperbolic 4-manifold without spin structures}

\begin{document}

\begin{abstract}
We build a non-compact, orientable, hyperbolic four-manifold of finite volume that does not admit any spin structure. 
\end{abstract}

\maketitle

\section*{Introduction}

It follows from a couple of works of Deligne and Sullivan \cite{DS,S} of the 1970s that every hyperbolic manifold $M$ is finitely covered by a stably parallelisable manifold $M'$. In particular, the Stiefel--Whitney classes satisfy $w_k(M') = 0$ for all $k > 0$. Unless otherwise stated, all manifolds in the paper are smooth, connected, and orientable (ie with $w_1 = 0$), and all hyperbolic manifolds are complete and of finite volume. 

The existence of hyperbolic $n$-manifolds that do not admit spin structures (ie with $w_2 \neq 0$) has been proved in 2020:  there are closed for all $n \geq 4$ \cite{MRS1} and cusped for all $n \geq 5$ \cite{LR}. Recall instead that surfaces are stably parallelisable and 3-manifolds are parallelisable.
Then, several examples of hyperbolic manifolds with non-trivial Stiefel--Whitney classes have been produced with different techniques \cite{C1,C2,LR,MRS2,RS1}, but the existence of cusped 4-manifolds with $w_2 \neq 0$ appears open. We fill here the gap:

\begin{thm} \label{thm:main}
    There exists a 
    cusped orientable (arithmetic) hyperbolic $4$-manifold $M$ 
    that does not admit any spin structure.
\end{thm}

Since $M$ is arithmetic and even-dimensional, by \cite[Lemma 4.2]{ERT} we can iteratively apply the embedding theorem of Kolpakov, Reid and Slavich \cite{KRS} as in \cite[Section 5]{MRS1}, to get a sequence of totally geodesic embeddings $M = \matH^4/\Gamma_4 \subset \matH^{5}/\Gamma_{5} \subset \ldots$ of $n$-manifolds with $\Gamma_n \subset \PSO(1,n;\matQ)$ commensurable with $\PO(1,n;\matZ)$. None of them admits a spin structure because an orientable hypersurface does not, so: 

\begin{cor} \label{cor}
    For every $n\ge4$, there exists a 
    cusped orientable (arithmetic) hyperbolic $n$-manifold 
    that does not admit any spin structure.
\end{cor}

This has already been proved by Long and Reid for $n \geq 5$ \cite{LR} as follows: (1) there is a closed flat 4-manifold $F^4$ with $w_2(F^4) \neq 0$ \cite{LRT, PS}, so $F^{n-1} = F^4 \times S^1 \times \ldots \times S^1$ has $w_2(F^{n-1}) \neq 0$ for all $n \geq 5$; (2) as every closed flat manifold, $F^{n-1}$ is diffeomorphic to a cusp section 
of a cusped hyperbolic manifold $M^n$ \cite{LR1,M1}, so as before $w_1(F^{n-1}) = 0,\ w_2(F^{n-1}) \neq 0 \implies w_2(M^n) \neq 0$. 

To prove Theorem \ref{thm:main}, we instead proceed as done in the closed case by Martelli, Slavich and the first author in \cite{MRS1} (see also \cite{MRS2}), explicitly constructing a hyperbolic 4-manifold $M$ satisfying a stronger condition: its intersection form is \emph{odd}; equivalently, there is a closed oriented surface $S \subset M$ with odd self-intersection $S \cdot S$ (the Euler number of the normal bundle). Then $w_2(M) \neq 0$ because the result of clashing $w_2(M)$ with the $\matZ/2\matZ$-homology class of $S$ is $S \cdot S \mod2$. Note that $S$ must necessarily be closed, 
otherwise $S \cdot S = 0$. Moreover, $S$ is not homologous to any immersed totally geodesic surface in $M$, since such surfaces have even self-intersection (see \cite{MRS1}).

As in \cite{MRS1,MRS2}, we build $M$ by gluing some copies of a right-angled hyperbolic polytope $P^4$ in such a way that $S$ is contained in the 2-skeleton of the tessellation. For this purpose, we need that $P^4$ has a compact 2-face $P^2$.
The only unbounded, right-angled, hyperbolic 4-polytope of finite volume with a compact 2-face that we know is introduced in Section \ref{sec:P}. It belongs to a continuous family of hyperbolic 4-polytopes discovered in 2010 by Kerckhoff--Storm \cite{KS}, further studied in \cite{MR} and later used for different purposes \cite{R,RS2,RS2',RS3}.
The polytope $P^4$ has 22 facets and octahedral symmetry. Its reflection group is arithmetic, and like for the well-known ideal 24-cell, is commensurable with the integral lattice $\PO(1,4;\matZ)$. {The 
manifold
$M$ 
belongs to this commensurability class.} 
We thank Leone Slavich for pointing out that a conjugate of $\Gamma_4 \cong \pi_1(M)$ lies in $\PSO(1,4;\matQ)$, which gives Corollary \ref{cor}.

Like in \cite{M,MRS1,MRS2}, we use some right-angled polytopes $P^2 \subset P^3 \subset P^4$ (where $P^n$ is a facet of $P^{n+1}$) to build some auxiliary hyperbolic manifolds with right-angled corners of increasing dimension. These objects have been fruitfully used in four- and five- dimensional hyperbolic geometry in the very last years \cite{BFS,C1,C3,R}. The surface $S$ is piecewise geodesic and tessellated by copies of $P^2$, and the cells of $M$ intersecting $S$ form a 4-manifold with right-angled corners $X$ (see Figure \ref{fig:N,X}--right). Like in \cite{MRS1,MRS2}, the construction ensures 
the following:

\begin{thm} \label{thm:geom-fin}
    There exists a geometrically finite hyperbolic $4$-manifold (of infinite volume) that covers a cusped hyperbolic manifold (of finite volume) and deformation retracts onto a closed surface with non-trivial normal bundle.
\end{thm}

Theorem \ref{thm:geom-fin} follows from the fact that $M$ contains $X$ as a convex submanifold, and the latter deformation retracts onto $S$. So $\pi_1(S)$ injects in $\pi_1(M)$ and induces a covering $\hat{M} \to M$ such that $\hat{M}$ is geometrically finite and diffeomorphic to the interior of $X$. For a proof, substitute ``compact'' with ``complete and finite-volume'' and ``convex cocompact'' with ``geometrically finite'' in the proof of \cite[Proposition 6, Corollary 8]{MRS2}.

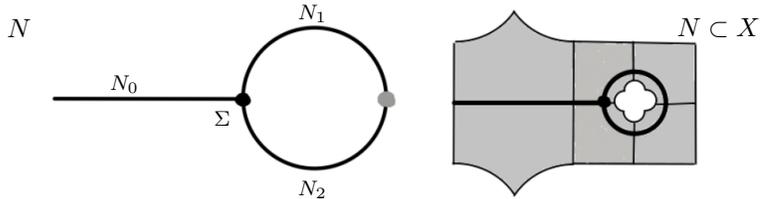
\begin{figure}
    \centering
    \begin{tikzpicture}
        \node[anchor=south west,inner sep=0] at (0,0) {\includegraphics[scale=2]{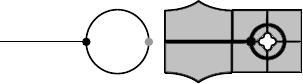}};
        \node at (0,2.85){$N$};
        \node at (9.3,2.85){$N \subset X$};
        \node at (1.4,1.7){\footnotesize $N_0$};
        \node at (4,2.75){\footnotesize $N_1$};
        \node at (4,.6){\footnotesize $N_2$};
        \node at (2.7,1.7){\footnotesize $\Sigma$};
    \end{tikzpicture}
 
    \caption{\footnotesize On the left, a schematic picture of the three-dimensional thickening $N = N_0 \cup N_1 \cup N_2$ of the piecewise geodesic surface $S = S_0 \cup S_1 \cup S_2$, where $S_i \subset N_i$ are totally geodesic manifolds with corners. It is not a manifold near the auxiliary surface with corners $\Sigma = N_0 \cap N_1 \cap N_2$ (represented by a black dot). On the right, the thickening $X$ of $N$: a 4-manifold with corners, neighbourhood of $S$ in $M$, tessellated by some copies of $P^4$ (represented by 10 gray pentagons)} \label{fig:N,X}
\end{figure}

The paper is organised as follows: the proof of Theorem \ref{thm:main} is summarised in Section \ref{sec:summary}, the polytope is introduced in Section \ref{sec:P}, and the construction is performed in Section \ref{sec:constr}.

\section{Summary} \label{sec:summary}

As already explained, like in \cite{MRS1,MRS2} for the compact case, our goal is to prove the following:

\begin{thm} \label{thm}
    There exists a 
    cusped, oriented, arithmetic, hyperbolic $4$-manifold $M$ that contains an oriented 
    surface $S$ with self-intersection $S \cdot S = 1$.
\end{thm}

%
%
A \emph{(hyperbolic) manifold with (right-angled) corners} is a complete hyperbolic manifold with boundary $X$, locally modelled on an orthant of $\matH^n$. The connected submanifolds with boundary that naturally stratify $\partial X$ are called \emph{faces}. We call \emph{facets} and \emph{corners} the $(n-1)$-dimensional and $(n-2)$-dimensional faces, respectively. Each face is naturally the image under a local isometry of a manifold with corners. These local isometries are all embeddings precisely when every corner is the intersection of two facets.

An $n$-manifold with corners and embedded facets $X$ is contained in a hyperbolic $n$-manifold $M$ without boundary, obtained in a standard way by iteratively doubling and re-doubling $X$ along its facets (see Section \ref{sec:M}). So, to prove Theorem \ref{thm}, we are reduced to building a cusped 4-manifold with corners $X$ with embedded faces and a surface $S \subset X$ such that $S \cdot S = 1$. 

The surface $S$ cannot be contained in an orientable 3-manifold in $M$, otherwise $S \cdot S = 0$. Similarly to \cite{MRS1}, it will instead be contained in a  ``locally Y-shaped piece'' $N$ obtained by gluing three 3-manifolds with corners $N_0$, $N_1$ and $N_2$ (also) along an isometric facet $\Sigma$ (see Figure \ref{fig:N,X}--left). The intersection $\Theta = \Sigma \cap S = \gamma_0 \cup \gamma_1 \cup \gamma_2$ is a theta-graph that trisects $S$ in three pieces $S_0$, $S_1$ and $S_2$, with $S_i$ properly embedded in $N_i$ and $\gamma_i = \Sigma \cap S_i$ a boundary component of $S_i$ (see Figure \ref{fig:sigma}). The 4-manifold with corners $X$ will be a thickening of $N$ (see Figure \ref{fig:N,X}--right), and will contain $S$ with $S \cdot S = \pm 1$ 
by construction (see Figure \ref{fig:isotopy}).

All $\Theta$, $\Sigma$, $S$ and $N$ will be contained in the skeleta of the tessellation of $X$ in copies of $P^4$. The auxiliary surface $\Sigma$ is totally geodesic, while $S$ is pleated. Moreover, $\Sigma$ and $S$ are tessellated by $P^2$'s and $N$ by $P^3$'s. Each $N_i$ is totally geodesic in $X$, and $N_0 \perp N_1, N_2$. The thickenings 
$S \subset N \subset X$ are built via the sequence $P^2 \subset P^3 \subset P^4$.

%

\section{The polytope} \label{sec:P}

We introduce here Kerckhoff and Storm's right-angled hyperbolic 4-polytope $P^4$ \cite{KS}. Let us identify the hyperbolic 4-space $\matH^4$ with the upper sheet of the hyperboloid $\langle x, x \rangle = -1$ in the Minkowski 5-space $\matR^{1,4}$. Here $\langle x, y \rangle = -x_0 y_0 + x_1 y_1 + \ldots + x_4 y_4$ for $x = (x_0, \ldots, x_4), y = (y_0, \ldots, y_4) \in \matR^{1,4}$. Given a spacelike vector $v \in \matR^{1,4}$, the inequality $\langle x, v \rangle \leq 0$ defines a half-space of $\matH^4$.
Let\footnote{In \cite{KS,MR}, $P^4$ is denoted by $P_t$, where $t = t_4 = \overline{t} = \sqrt3/3$.} $P^4 \subset \matH^4$ be the intersection of the 22 half-spaces given by the vectors in Table \ref{table}. It is an unbounded, right-angled polytope of finite volume \cite[Proposition 13.1]{KS}.

\begin{table}  
\begin{eqnarray*}
E_1  = \big( \sqrt{2},+1,+1,+1,+\sqrt3 \big) \quad & H_1  = \big( \sqrt{2},-1,-1,-1,+\sqrt3/3 \big) \quad & C_{12} = \big( 1,+\sqrt{2},0,0,0 \big)\\
E_1' = \big( \sqrt{2},-1,-1,-1,-\sqrt3 \big) \quad & H_1' = \big( \sqrt{2},+1,+1,+1,-\sqrt3/3 \big) \quad & C_{34} = \big( 1,-\sqrt{2},0,0,0 \big)\\
E_2  = \big( \sqrt{2},+1,-1,-1,+\sqrt3 \big) \quad & H_2  = \big( \sqrt{2},-1,+1,+1,+\sqrt3/3 \big) \quad & C_{13} = \big( 1,0,+\sqrt{2},0,0 \big) \\
E_2' = \big( \sqrt{2},-1,+1,+1,-\sqrt3 \big) \quad & H_2' = \big( \sqrt{2},+1,-1,-1,-\sqrt3/3 \big) \quad & C_{24} = \big( 1,0,-\sqrt{2},0,0 \big)\\
E_3  = \big( \sqrt{2},-1,+1,-1,+\sqrt3 \big) \quad & H_3  = \big( \sqrt{2},+1,-1,+1,+\sqrt3/3 \big) \quad & C_{14} = \big( 1,0,0,+\sqrt{2},0 \big)\\
E_3' = \big( \sqrt{2},+1,-1,+1,-\sqrt3 \big) \quad & H_3' = \big( \sqrt{2},-1,+1,-1,-\sqrt3/3 \big) \quad & C_{23} = \big( 1,0,0,-\sqrt{2},0 \big)\\
E_4  = \big( \sqrt{2},-1,-1,+1,+\sqrt3 \big) \quad & H_4  = \big( \sqrt{2},+1,+1,-1,+\sqrt3/3 \big) \quad & \\
E_4' = \big( \sqrt{2},+1,+1,-1,-\sqrt3 \big) \quad & H_4' = \big( \sqrt{2},-1,-1,+1,-\sqrt3/3 \big) \quad & 
\end{eqnarray*}
\caption{\footnotesize The spacelike vectors of $\matR^{1,4}$ that define the polytope $P^4 \subset \matH^4$.}\label{table}
\end{table}

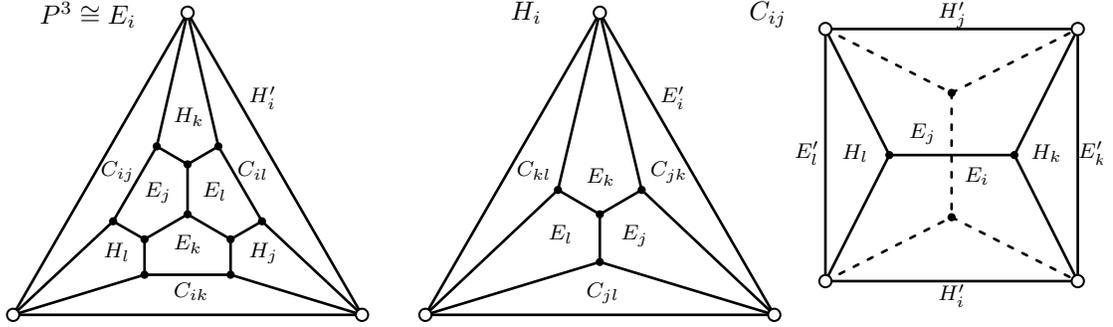
\begin{figure}
\begin{tikzpicture}[scale=.773,line cap=round,line join=round,>=triangle 45,x=1.0cm,y=1.0cm]
\clip(-0.18,-0.21) rectangle (18.74,5.49);
\draw [line width=1pt] (0,0)-- (1.72,1.61);
\draw [line width=1pt] (1.72,1.61)-- (2.26,1.3);
\draw [line width=1pt] (0,0)-- (2.26,0.69);
\draw [line width=1pt] (2.26,0.69)-- (2.26,1.3);
\draw [line width=1pt] (3,5.2)-- (2.47,2.9);
\draw [line width=1pt] (2.47,2.9)-- (3,2.59);
\draw [line width=1pt] (3,5.2)-- (3.53,2.9);
\draw [line width=1pt] (3.53,2.9)-- (3,2.59);
\draw [line width=1pt] (6,0)-- (4.28,1.61);
\draw [line width=1pt] (4.28,1.61)-- (3.74,1.3);
\draw [line width=1pt] (3.74,0.69)-- (3.74,1.3);
\draw [line width=1pt] (6,0)-- (3.74,0.69);
\draw [line width=1pt] (1.72,1.61)-- (2.47,2.9);
\draw [line width=1pt] (3.53,2.9)-- (4.28,1.61);
\draw [line width=1pt] (3.74,0.69)-- (2.26,0.69);
\draw [line width=1pt] (2.26,1.3)-- (3,1.73);
\draw [line width=1pt] (3,1.73)-- (3.74,1.3);
\draw [line width=1pt] (3,1.73)-- (3,2.59);
\draw [line width=1pt] (10.09,5.2)-- (10.81,2.15);
\draw [line width=1pt] (10.81,2.15)-- (13.09,0);
\draw [line width=1pt] (10.09,5.2)-- (9.37,2.15);
\draw [line width=1pt] (7.09,0)-- (10.09,0.91);
\draw [line width=1pt] (10.09,0.91)-- (13.09,0);
\draw [line width=1pt] (7.09,0)-- (9.37,2.15);
\draw [line width=1pt] (9.37,2.15)-- (10.09,1.73);
\draw [line width=1pt] (10.09,1.73)-- (10.81,2.15);
\draw [line width=1pt] (10.09,1.73)-- (10.09,0.91);
\draw [line width=1pt] (7.09,0)-- (10.09,5.2);
\draw [line width=1pt] (10.09,5.2)-- (13.09,0);
\draw [line width=1pt] (13.09,0)-- (7.09,0);
\draw [line width=1pt] (6,0)-- (3,5.2);
\draw [line width=1pt] (3,5.2)-- (0,0);
\draw [line width=1pt] (0,0)-- (6,0);
\draw [line width=1pt] (13.97,4.92)-- (18.31,4.92);
\draw [line width=1pt] (18.31,4.92)-- (18.31,0.58);
\draw [line width=1pt] (18.31,0.58)-- (13.97,0.58);
\draw [line width=1pt] (13.97,0.58)-- (13.97,4.92);
\draw [line width=1pt,dash pattern=on 2pt off 4pt] (13.97,4.92)-- (16.14,3.82);
\draw [line width=1pt,dash pattern=on 2pt off 4pt] (16.14,3.82)-- (18.31,4.92);
\draw [line width=1pt,dash pattern=on 2pt off 4pt] (13.97,0.58)-- (16.14,1.68);
\draw [line width=1pt,dash pattern=on 2pt off 4pt] (16.14,1.68)-- (18.31,0.58);
\draw [line width=1pt,dash pattern=on 2pt off 4pt] (16.14,3.82)-- (16.14,1.68);
\draw [line width=1pt] (15.07,2.75)-- (17.22,2.75);
\draw [line width=1pt] (13.97,4.92)-- (15.07,2.75);
\draw [line width=1pt] (15.07,2.75)-- (13.97,0.58);
\draw [line width=1pt] (17.22,2.75)-- (18.31,4.92);
\draw [line width=1pt] (17.22,2.75)-- (18.31,0.58);
\draw (0.3,5.55) node[anchor=north west] {$P^3 \cong E_i$};
\draw (2.6,1.54) node[anchor=north west] {\footnotesize $E_k$};
\draw (2.1,2.44) node[anchor=north west] {\footnotesize $E_j$};
\draw (3.1,2.44) node[anchor=north west] {\footnotesize $E_l$};
\draw (2.6,0.72) node[anchor=north west] {\footnotesize $C_{ik}$};
\draw (3.7,2.8) node[anchor=north west] {\footnotesize $C_{il}$};
\draw (1.35,2.8) node[anchor=north west] {\footnotesize $C_{ij}$};
\draw (1.4,1.4) node[anchor=north west] {\footnotesize $H_l$};
\draw (3.9,1.4) node[anchor=north west] {\footnotesize $H_j$};
\draw (2.6,3.73) node[anchor=north west] {\footnotesize $H_k$};
\draw (3.9,4.1) node[anchor=north west] {\footnotesize $H'_i$};
\draw (8.4,5.55) node[anchor=north west] {$H_i$};
\draw (9.7,2.69) node[anchor=north west] {\footnotesize $E_k$};
\draw (9.05,1.7) node[anchor=north west] {\footnotesize $E_l$};
\draw (10.97,4.1) node[anchor=north west] {\footnotesize $E'_i$};
\draw (8.5,2.8) node[anchor=north west] {\footnotesize $C_{kl}$};
\draw (10.8,2.8) node[anchor=north west] {\footnotesize $C_{jk}$};
\draw (9.7,0.72) node[anchor=north west] {\footnotesize $C_{jl}$};
\draw (10.3,1.7) node[anchor=north west] {\footnotesize $E_j$};
\draw (12.5,5.55) node[anchor=north west] {$C_{ij}$};
\draw (16.2,2.7) node[anchor=north west] {\footnotesize $E_i$};
\draw (15.25,3.45) node[anchor=north west] {\footnotesize $E_j$};
\draw (17.35,3.1) node[anchor=north west] {\footnotesize $H_k$};
\draw (14.1,3.1) node[anchor=north west] {\footnotesize $H_l$};
\draw (15.75,5.55) node[anchor=north west] {\footnotesize $H'_j$};
\draw (13.28,3.15) node[anchor=north west] {\footnotesize $E'_l$};
\draw (18.1615,3.15) node[anchor=north west] {\footnotesize $E'_k$};
\draw (15.75,0.69) node[anchor=north west] {\footnotesize $H'_i$};
\begin{scriptsize}
\fill [color=black] (3,1.73) circle (2pt);
\fill [color=black] (1.72,1.61) circle (2pt);
\fill [color=black] (2.26,1.3) circle (2pt);
\fill [color=black] (0,0) circle (3.5pt);
\fill [color=white] (0,0) circle (2.5pt);
\fill [color=black] (2.26,0.69) circle (2pt);
\fill [color=black] (2.47,2.9) circle (2pt);
\fill [color=black] (3,5.2) circle (3.5pt);
\fill [color=white] (3,5.2) circle (2.5pt);
\fill [color=black] (3.53,2.9) circle (2pt);
\fill [color=black] (3,2.59) circle (2pt);
\fill [color=black] (4.28,1.61) circle (2pt);
\fill [color=black] (3.74,1.3) circle (2pt);
\fill [color=black] (6,0) circle (3.5pt);
\fill [color=white] (6,0) circle (2.5pt);
\fill [color=black] (3.74,0.69) circle (2pt);
\fill [color=black] (10.09,1.73) circle (2pt);
\fill [color=black] (10.81,2.15) circle (2pt);
\fill [color=black] (10.09,5.2) circle (3.5pt);
\fill [color=white] (10.09,5.2) circle (2.5pt);
\fill [color=black] (10.09,0.91) circle (2pt);
\fill [color=black] (13.09,0) circle (3.5pt);
\fill [color=white] (13.09,0) circle (2.5pt);
\fill [color=black] (7.09,0) circle (3.5pt);
\fill [color=white] (7.09,0) circle (2.5pt);
\fill [color=black] (9.37,2.15) circle (2pt);
\fill [color=black] (13.97,4.92) circle (3.5pt);
\fill [color=white] (13.97,4.92) circle (2.5pt);
\fill [color=black] (18.31,4.92) circle (3.5pt);
\fill [color=white] (18.31,4.92) circle (2.5pt);
\fill [color=black] (13.97,0.58) circle (3.5pt);
\fill [color=white] (13.97,0.58) circle (2.5pt);
\fill [color=black] (18.31,0.58) circle (3.5pt);
\fill [color=white] (18.31,0.58) circle (2.5pt);
\fill [color=black] (16.14,3.82) circle (2pt);
\fill [color=black] (16.14,1.68) circle (2pt);
\fill [color=black] (15.07,2.75) circle (2pt);
\fill [color=black] (17.22,2.75) circle (2pt);
\end{scriptsize}
\end{tikzpicture}
\caption{\footnotesize The extremal, half-height and central facets $E_i \cong P^3$, $H_i$ and $C_{ij}$ of $P^4$, where $\{ i,j,k,l \} = \{ 1,2,3,4 \}$. The ideal vertices are in white. Note the compact pentagon $E_i \cap E_j \cong P^2$.}\label{fig:P}
\end{figure}

Note that the isometry $a$ defined by $a(x_0, x_1, \ldots, x_4) = (x_0, -x_1, \ldots, -x_4)$ is a symmetry of $P^4$. Moreover, the notation (taken from \cite{R}) is such that $E_i' = a(E_i)$, $H_i' = a(H_i)$, and $C_{ij} = a(C_{kl})$ for all distinct $i,j,k,l$.
The combinatorics of $P^4$ has been studied in detail in \cite[Proposition 3.16]{MR}. Each vector in Table \ref{table} corresponds to a facet of $P^4$, denoted with the same symbol. The 22 facets, depicted in Figure \ref{fig:P}, are partitioned up to symmetry into three sets:\footnote{In \cite{KS,MR}, these are called: the 
``positive walls'', the 
``negative walls'', and the ``letter walls'', respectively.} 
\begin{enumerate}
    \item the \emph{extremal facets} $E_1, E_2, E_3, E_4, E'_1, E'_2, E'_3, E'_4$,
    \item the \emph{half-height facets} $H_1, H_2, H_3, H_4, H'_1, H'_2, H'_3, H'_4$,
    \item the \emph{central facets} $C_{12}, C_{13}, C_{14}, C_{23}, C_{24}, C_{34}$.
\end{enumerate}

\begin{lemma} \label{lem:aithm}
    Every combinatorial automorphism of $P^4$ is realised by an isometry of $P^4$, and every hyperbolic orbifold $O$ tessellated by finitely-many copies of $P^4$ is commensurable with $\matH^4/\PO(1,4;\matZ)$.
\end{lemma}
\begin{proof}
    The poof of the first statement (relying on \cite[Proposition 2.4]{RS3} and \cite[Lemma 4.15]{MR}) is the same of \cite[Lemma 1.2]{R} by \cite[Section 3.2 and Proposition 3.16]{MR}. In particular (see Figure \ref{fig:P}), every isometry between two facets of $P^4$ is the restriction of an isometry of $P^4$. Since, by hypothesis, $O$ can be obtained by gluing the facets of some copies $P^4$ in  pairs via isometries, $O$ covers the orbifold $P^4/\Isom(P^4)$, and so it is commensurable with $P^4 = \matH^4 / \Gamma$. The reflection group $\Gamma < \PO(1,n)$ of $P^4$ is arithmetic \cite[Theorem 13.2]{KS} and commensurable with 
    $\PO(1,4;\matZ)$ \cite[Proposition 4.25]{MR}.
\end{proof}

Note from Figure \ref{fig:P} that the compact 2-faces of $P^4$ are 12 isometric pentagons $E_i \cap E_j$, $E_i' \cap E_j'$, $i \neq j$. Defining 
$$P^2 = E_1 \cap E_2\ \mbox{and}\ P^3 = E_1,$$
we have a sequence of right-angled polytopes:
$$P^2 \subset P^3 \subset P^4.$$
We shall think of $P^{n+1}$ as sitting above its \emph{bottom facet} $P^n$, and call the remaining facets \emph{vertical facets} and \emph{top facets}, depending on whether they are adjacent to $P^n$ or not, respectively. For example, $P^3$ has 5 vertical facets and 4 top facets, while $P^4$ has 10 vertical facets and 11 top facets.

\section{The construction} \label{sec:constr}

In this section, we prove Theorem \ref{thm}. We first build the auxiliary surface with corners $\Sigma$, and thicken it to a 3-manifold with corners $\Sigma^\thick$ homeomorphic to $\Sigma \times [0,1]$. Then, we build the 3-manifolds with corners $N_0$, $N_1$ and $N_2$ by gluing some of the top facets of $\Sigma^\thick$ in three different ways, and the 3-manifold with corners $N_{12}$ by gluing together $N_1$ and $N_2$. After that, we glue the 3-manifolds with corners $N_0$ and $N_{12}$ and thicken the resulting ``locally Y-shaped piece'' $N$ to a 4-manifold with corners $X$. Then we study $X$, and finally build the 4-manifold $M$.

\subsection{The surface with corners $\Sigma$ and its thickening $\Sigma^\thick$} \label{sec:sigma}

\begin{figure}
    \centering
    \begin{tikzpicture}
        \node[anchor=south west,inner sep=0] at (0,0) {\includegraphics[scale=1.2]{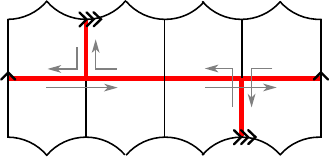}};
        \node at (-0.5,2.7){$\Sigma$};
        \node at (2.3,1.1){\footnotesize $\gamma_0$};
        \node at (5.5,1.1){\footnotesize $\gamma_0$};
        \node at (1.3,2){\footnotesize $\gamma_2$};
        \node at (2.3,2){\footnotesize $\gamma_1$};
        \node at (5.4,2){\footnotesize $\gamma_2$};
        \node at (4.6,2){\footnotesize $\gamma_1$};
    \end{tikzpicture}
    
    \caption{\footnotesize The surface $\Sigma$ with corners obtained by gluing 8 copies of the right-angled pentagon $P^2$ (four edges of the big dodecagon are glued in pairs as indicated by the black arrows). It is a holed torus, and deformation retracts onto the red theta-graph $\Theta = \gamma_0 \cup \gamma_1 \cup \gamma_2$. The three red oriented curves $\gamma_0, \gamma_1$ and $\gamma_2$ run as indicated by the gray arrows.}
    \label{fig:sigma}
\end{figure}
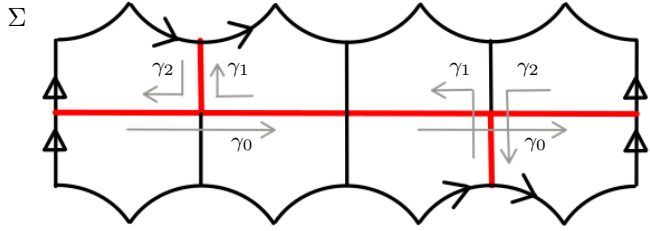

Let $\Sigma$ be the surface with corners obtained by gluing in pairs some edges of $8$ copies of $P^2$ 
as indicated in Figure \ref{fig:sigma}. Topologically, $\Sigma$ is a once-holed torus. Consider the three oriented curves $\gamma_0$, $\gamma_1$ and $\gamma_2$ in the 1-skeleton of $\Sigma$ as in Figure \ref{fig:sigma}. The surface $\Sigma$ is a thickening of the theta-graph $\Theta = \gamma_0 \cup \gamma_1 \cup \gamma_2$.

We now place a copy of $P^3$ ``above'' each $P^2$ in $\Sigma$, to get a 3-manifold with corners $\Sigma^\thick$ homeomorphic to $\Sigma \times [0,1]$: the vertical faces of the $P^3$'s containing the paired edges of the $P^2$'s in $\Sigma$ are glued correspondingly
. So $\Sigma^\thick$ has three types of facets: the \emph{bottom facet} $\Sigma$, and the \emph{vertical} and \emph{top facets} tessellated by the facets of $P^3$ of the corresponding type.

\begin{figure}
    \centering
    \includegraphics[scale=1.2]{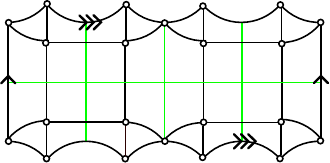}
    \caption{\footnotesize The top of the $3$-manifold with corners $\Sigma^\thick$
    . The green lines indicate its tessellation into 8 copies of $P^3$. As usual, the ideal vertices are in white.}
    \label{fig:sigma-thick}
\end{figure}

Recall that $P^3$ has four top facets: an ideal triangle, a quadrilateral with one ideal vertex, and two quadrilaterals with two consecutive ideal vertices. The top facets of $\Sigma^\thick$, tessellated by the above polygons, are $8$ ideal triangles, $4$ ideal rectangles and $3$ ideal hexagons, pleated with right angles along the pattern showed in Figure \ref{fig:sigma-thick}. 

\subsection{The 3-manifolds with corners $N_0$, $N_1$, $N_2$ and $N_{12}$} \label{sec:Ni}

\begin{figure}
    \centering
    \includegraphics[scale=1.2]{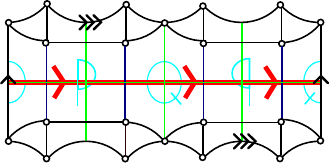}
    \caption{\footnotesize The 3-manifold with corners $N_0$ is built by gluing some top facets of $\Sigma^\thick$ as indicated by the blue letters P and Q. It has 5 top facets. The four vertical blue edges are glued making an angle of $2\pi$. }
    \label{fig:N0}
\end{figure}

\begin{figure}
    \centering
    \includegraphics[scale=1.2]{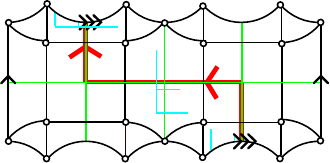}
    \caption{\footnotesize The top of the 3-manifold with corners $N_1$.}
    \label{fig:N1}
\end{figure}

\begin{figure}
    \centering
    \includegraphics[scale=1.2]{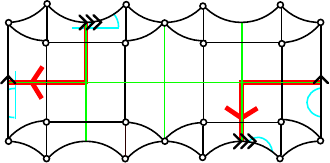}
    \caption{\footnotesize The top of the 3-manifold with corners $N_2$.}
    \label{fig:N2}
\end{figure}

Let $N_0$, $N_1$ and $N_2$ be obtained by gluing some top facets of $\Sigma^\thick$ in pairs 
as indicated by Figures \ref{fig:N0}, \ref{fig:N1} and \ref{fig:N2}, respectively. Specifically, for $N_0$ we have paired the two hexagons with label Q and two quadrilaterals with label P, for $N_1$ the two hexagons with label F, and for $N_2$ the two hexagons with label P.

Figure \ref{fig:N0} helps to verify that $N_0$ is a 3-manifold with corners and embedded facets: the four glued corners are cyclically glued together in the interior of $N_0$, and each of the remaining corners is right-angled and belongs to two distinct facets. Moreover, the 8 copies of $P^3$ in $N_0$ that are adjacent to $\Sigma$ are distinct. The check for $N_1$ and $N_2$ is even simpler, and is left to the reader. 

For $i=0,1,2$, consider the surface with corners $S_i' \cong \gamma_i \times [0,1]$ in the 2-skeleton of $\Sigma^\thick$, tessellated by the vertical pentagons that have an edge in $\gamma_i \subset \Sigma \subset \Sigma^\thick$. The red line in Figures \ref{fig:N0}, \ref{fig:N1} and \ref{fig:N2} is the top of $S_i'$. We call $S_i$ the surface in $N_i$ obtained from $S_i'$ after the gluing. Both $N_i$ and $S_i$ are orientable, since the gluings reverse the orientation of both the glued polygons and the red curve.

\begin{figure}
    \centering
    \includegraphics[scale=1.2]{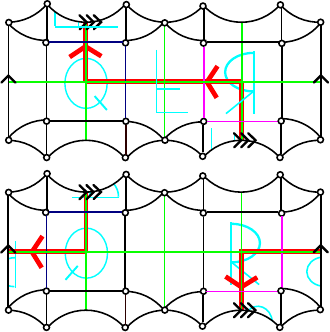}
    \caption{\footnotesize The top of the 3-manifold with corners $N_{12}$, obtained by pairing some top facets of $N_1$ (top) and $N_2$ (bottom) as indicated by the symbols P, Q, R and F, and the two bottom facets. It has 15 top facets. The four blue (resp. pink) edges are glued making an angle of $2\pi$.}
    \label{fig:N12}
\end{figure}

We conclude by gluing together $N_1$ and $N_2$ as follows: we glue their two bottom facets (copies of $\Sigma$) 
and some of their top facets as in Figure \ref{fig:N12}. In particular, we have paired the two top quadrilaterals with label Q and the ones with label R (whereas the hexagons with labels P and F have already been paired before while building $N_1$ and $N_2$). We call $N_{12}$ the resulting 3-manifold with corners. 
Again, it is easy to check that $N_{12}$ is an orientable 3-manifold with corners and embedded facets, that the 16 copies of $P^3$ in $N_{12}$ incident to $\Sigma \subset N_{12}$ are distinct, and that $S_{12} = S_1 \cup S_2$ is an orientable surface embedded in $N_{12}$ with $\partial S_{12} = \gamma_1 \sqcup \gamma_2$.

\subsection{The spine $N$ and its thickening $X$} \label{sec:N}

Let $N$ be obtained by gluing $N_0$ and $N_{12}$ 
along their two isometric copies of $\Sigma$: the bottom facet of $N_0$ and the properly embedded surface in $N_{12}$ obtained by identifying the two bottom facets of $N_1$ and $N_2$. It is not a manifold (see Figure \ref{fig:N,X}--left). 

We now want to thicken $N$ to a 4-manifold with corners $X$ in which $N_0$ and $N_{12}$ are totally geodesic and orthogonal. Similarly to \cite{MRS1,MRS2}, this can be done in two steps.

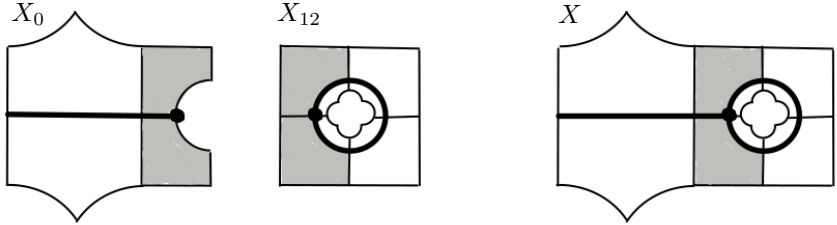
\begin{figure}
    \centering
    \begin{tikzpicture}
        \node[anchor=south west,inner sep=0] at (0,0) {\includegraphics[scale=1.5]{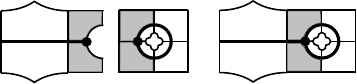}};
        \node at (1.3,2.4){$X_0$};
        \node at (4,2.4){$X_{12}$};
        \node at (7.3,2.4){$X$};
    \end{tikzpicture}
    \caption{\footnotesize A schematic picture of $X$ (right), obtained by gluing the thickenings $X_0$ and $X_{12}$ of $N_0$ and $N_{12}$ (left). The pentagons, thick segment, circle and dot represent the copies of $P^4$, the 3-manifolds $N_0$ and $N_{12}$, and the surface $\Sigma$, respectively.}
    \label{fig:X}
\end{figure}

We first thicken $N_0$ and $N_{12}$ separately: we place two copies of $P^4$ on every copy of $P^3$, one ``below" and the other ``above", and get two 4-manifolds with corners $X_0$ and $X_{12}$; see Figure \ref{fig:X}--left. These are obtained by pairing some vertical facets of some copies of $P^4$
.

Then, we identify in pairs the copies of $P^4$ in $X_0$ incident to $\Sigma$ with the copies of $P^4$ in $X_{12}$ incident to $\Sigma$ from below, as in Figure \ref{fig:X}--right. We can do this since for every pentagonal face $F$ of $P^4$, there exists an isometry of $P^4$ that exchanges the two facets (isometric to $P^3$) that share $F$. The resulting complex $X$ contains $N$ as desired.

\subsection{The 4-manifold with corners $X$} \label{sec:X}

By construction, $X$ is a complete and orientable hyperbolic $4$-manifold with boundary. Since it is tessellated by copies of $P^4$, which is right angled, a priori the angles at the corners are multiples of $\pi/2$.
Our aim is now to show that the angles are $\pi/2$ and the facets are embedded. 

\begin{prop} 
    The thickening $X$ is a hyperbolic manifold with right-angled corners.  
\end{prop}
\begin{proof}
Every copy of a pentagonal face $F_i \cap F_j \cong P^2$ of $P^4$ contained in $\partial X$ is shared by 
one or two copies of $P^4$ in $X$. Indeed, as we see from Figure \ref{fig:P}, if $F$ and $F'$ are two facets of $P^4$ such that $F \cap E_i \ne \emptyset$, $F \cap E_j = \emptyset$, $F' \cap E_i = \emptyset$, $F' \cap E_j \ne \emptyset$, then $F \cap F' = \emptyset$. Therefore, such a copy of $P^2$ is either contained in a corner with angle $\pi/2$, or in a facet (with angle $\pi$).
\end{proof}

We now want to show that $X$ has embedded facets. We begin by showing that $X_0$ and $X_{12}$ have embedded facets. Let $Y$ be a facet of $X_0$ or $X_{12}$. 
It is a union of copies of a facet $F$ of $P^4$. Consider a corner $C$ of $X$ contained in $Y$. It is not possible that both sides of $C$ are in $Y$. Indeed, $P^4$ is right angled, hence both sides of $C$ are in the same copy of $P^4$ and, of course, there is only one facet $F$ in $P^4$.

\begin{prop}
    The facets of $X$ are embedded.
\end{prop}
\begin{proof}
Let $i \colon X_0 \to X$ and $j \colon X_{12} \to X$ be the natural inclusion embeddings. Let $Y$ be a facet of $X$. If $Y$ is entirely contained in $i(X_0)$ or $j(X_{12})$, then we easily conclude as 
above for $X_0$ and $X_{12}$. Otherwise, $Y \cap i(X_0)$ is union of copies of the facet $F$ of $P^4$ and $Y \cap j(X_{12})$ is union of copies of the facet $f(F)$ of $P^4$, where $f$ is the isometry used for the identification in the construction of $X$ starting from $X_0$ and $X_{12}$. We conclude as 
above for $X_0$ and $X_{12}$, since $P^4$ has only one facet $F$ and one facet $f(F)$.
\end{proof}

The construction ensures the following.

\begin{prop}
    The surface $S = S_0 \cup S_{12}$ has self-intersection $\pm 1$ in $X$.
\end{prop}
\begin{proof}
    We isotope $N$ inside a regular neighbourhood $U$ of $N$ in $X$ as follows. Say that $U = U_0 \cup U_{12}$ for two tubular neighbourhoods $U_0$ and $U_{12}$ of $N_0$ and $N_{12}$. The latter are two-sided. Call $U_+$ and $U_-$ the two sides of $U_{12}$, with $U_+ \cup U_- = U_{12}$ and $U_+ \cap U_- = N_{12}$.
    
    We first move $N$ in one direction as in Figure \ref{fig:isotopy}--left, obtaining an isotopic copy $N' = N_0' \cup N_{12}'$ of $N$ transverse to it, with $N' \cap N \subset U_{12}$. Then, to remove the intersection with $N_{12}'$, we ``push'' $N' \cap U_+$ in the interior of $U_-$ as in Figure \ref{fig:isotopy}--right, obtaining $N'' = N_0' \cup N_{12}''$. 

    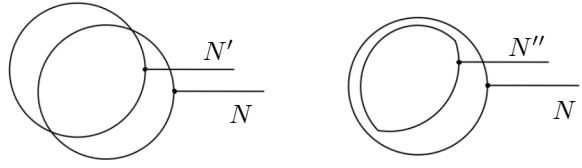
\begin{figure}
    \centering
    \begin{tikzpicture}
        \node[anchor=south west,inner sep=0] at (0,0) {\includegraphics[scale=1.4]{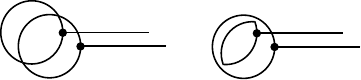}};
        \node at (4,1){$N$};
        \node at (3.6,1.4){$N'$};
        \node at (8.6,1){$N$};
        \node at (8.3,1.4){$N''$};
    \end{tikzpicture}
    \caption{\footnotesize On the left, $N$ and its isotopic copy $N'$. On the right, $N$ and its isotopic copy $N''$, which transversely intersects $N$ in the point $S \cap S''$, so $S \cdot S = S \cdot S'' = \pm 1$.}
    \label{fig:isotopy}
\end{figure}

    Then the surface $\Sigma''' = N \cap N'' = N_{12} \cap N_0'$ is a surface parallel to $\Sigma$ and $\Sigma''$. Moreover, $S$ and its isotopic copy $S'' \subset N''$ intersect transversely at one point, corresponding to the transverse intersection of the simple closed curves $\gamma_0''' = S_0 \cap \Sigma'''$ and $\gamma_2''' = S_2 \cap \Sigma'''$ in $\Sigma'''$.  Therefore $S \cdot S = S \cdot S'' = \pm 1$.
\end{proof}

\subsection{The 4-manifold $M$} \label{sec:M}

Let $Y_1, \ldots, Y_m$ the facets of $X$. We now double $X$ along $Y_1$, then double the result along the copies of $Y_2$, and continue iteratively, until we get a $4$-manifold $M$ without boundary tessellated by $2^m$ copies of $X$.
Since the facets of $X$ are embedded, $M$ is hyperbolic manifold (see eg \cite[Proposition 6]{M}). 
It is orientable because $X$ is orientable, and arithmetic by Lemma \ref{lem:aithm}. To complete the proof of Theorem \ref{thm}, it suffices to choose the orientation of $M$ such that $S \cdot S = +1$.

\end{document}